\documentclass[11pt]{amsart}
\font\emailfont=cmtt10

\usepackage[left=1.5 in, right=1.5 in, top=1 in, bottom=1 in]{geometry}
\usepackage{amsmath,amsthm,amsfonts,amscd,flafter,epsf}
\usepackage{graphics,todonotes,appendix}
\usepackage{epsfig}
\usepackage{psfrag}
\usepackage[all,cmtip]{xy}

\newcommand\commentable[1]{#1}

\newtheorem{thm}{Theorem}

\newtheorem*{Adjunction Inequality}{Adjunction Inequality}

\newtheorem{theorem}{Theorem}[section]
\newtheorem{prop}[theorem]{Proposition}
\newtheorem{corollary}[thm]{Corollary}

\newtheorem{conj}[thm]{Conjecture}

\theoremstyle{definition}
\newtheorem{example}[theorem]{Example}
\newtheorem{defn}[theorem]{Definition}
\newtheorem{remark}[theorem]{Remark}

\def\endproofof{\relax\ifmmode\expandafter\endproofmath\else
  \unskip\nobreak\hfil\penalty50\hskip.75em\hbox{}\nobreak\hfil\bull
  {\parfillskip=0pt \finalhyphendemerits=0 \bigbreak}\fi}
\def\endproofofmath$${\eqno\bull$$\bigbreak}

\def\endproof{\relax\ifmmode\expandafter\endproofmath\else
  \unskip\nobreak\hfil\penalty50\hskip.75em\hbox{}\nobreak\hfil\bull
  {\parfillskip=0pt \finalhyphendemerits=0 \bigbreak}\fi}
\def\endproofmath$${\eqno\bull$$\bigbreak}
\def\bull{\vbox{\hrule\hbox{\vrule\kern3pt\vbox{\kern6pt}\kern3pt\vrule}\hrule}}

\newcommand{\R}{\mathbb{R}}

\newcommand{\Z}{\mathbb{Z}}

\newcommand\SpinC{\mathrm{Spin}^c}

\newcommand\relspinc{\underline{\spinc}}

\newcommand\ModSphere{\ModFlow\left({\mathbb S}\longrightarrow 
\Sym^{g-1}(\Sigma_{1})\times \Sym^2(\Sigma_{2})\right)}
\newcommand\ModSpheres\ModSphere

\newcommand\CFa{\widehat{CF}}

\newcommand\HFa{\smash{\widehat{HF}}}

\newcommand\UnparModSp{\widehat \ModSp}
\newcommand\UnparModFlow\UnparModSp
\newcommand\Mod\ModSp

\newcommand{\spinc}{\mathfrak s}

\newcommand{\spinct}{\mathfrak t}

\newcommand\ModMaps{\mathcal M}
\newcommand\ModSp\ModMaps

\newcommand\spincrel\relspinc

\newcommand\CFK{CFK}

\newcommand\CFKinf{\CFK^{\infty}}

\newcommand\Dual{\mathcal D}
\newcommand\Duality\Dual

\newcommand\os{{Ozsv{\'a}th-Szab{\'o}}}

\newcommand\Sym{\mathrm{Sym}}

\newcommand\tb{\operatorname{tb}}
\newcommand\rot{\operatorname{rot}}

\newcommand\selflink{\operatorname{sl}}

\newcommand\lk{\operatorname{lk}}

\newcommand\image{\operatorname{Im}}

\commentable{

\title[{Relative Thom conjectures, symplectic and beyond}] 
{Relative Thom conjectures, symplectic and beyond}

}

\author[Matthew Hedden]{Matthew Hedden}
\address{Department of
Mathematics, Michigan State University, MI \newline
\indent{\emailfont{mhedden@math.msu.edu}}}

\author{Katherine Raoux}
\address{Department of
Mathematical Sciences, University of Arkansas, AR \newline
\indent{\emailfont{kraoux@uark.edu}}}

\begin{document}

\begin{abstract}
    We establish a criterion that ensures a bounded almost complex curve in a bounded almost complex 4-manifold minimizes genus amongst all smooth surfaces that share its homology class and the transverse link on its boundary.  An immediate corollary affirms the relative symplectic Thom conjecture and, moreover, yields obstructions coming from knot Floer homology to a link bounding a symplectic surface in a symplectic filling.  Our results are applicable to knots in  manifolds equipped with plane fields that admit no symplectic fillings; for instance, we show that symplectic surfaces in a thickening of any contact 3-manifold with non-zero \os\ invariant minimize slice genus for their boundary.  We conjecture that this phenomenon occurs precisely when the contact structure is tight, which would imply that tightness can be viewed as a symplecto-geometric notion.
\end{abstract}
\maketitle

\section{Introduction}
The Thom conjecture asserts a surprising aspect of smooth orientable surfaces in the complex projective plane: geometric complexity (holomorphicity) implies topological simplicity (genus minimality). First established by Kronheimer and Mrowka \cite{KMThom} (cf. \cite{MST}), this principle has become a paradigm in low-dimensional topology, culminating in Ozsv\'ath and Szab\'o's

\begin{thm} [Symplectic Thom Conjecture \cite{SympThom}] \label{thm:SympThom} If $\Sigma$ is a smoothly embedded symplectic surface in a closed symplectic $4$-manifold  $(W,\omega)$,  then \[-\chi(\Sigma)\le -\chi(\Sigma')\] for any smoothly embedded $\Sigma'\subset W$ with $[\Sigma]=[\Sigma']\in H_2(W)$.
\end{thm}

The purpose of this article is to push Thom-type conjectures beyond the holomorphic and symplectic realms. Our main result concerns  almost complex cobordisms $(W,J)$ to a 3-manifold $Y$. The almost complex structure induces a plane field $P_J=TY\cap JTY$ on $Y$, and an associated map on Floer cohomology \[F_{W,\kappa_J}^*:\HFa^*(Y)\to \HFa^*(Y_{in})\] where $Y_{in}$ is the (possibly empty) incoming end of $W$ and $\kappa_J$ is the $\SpinC$ structure associated to $J$. Given this, we show:
\begin{thm}[Relative Almost Complex Thom Conjecture] \label{thm:genRelativeSymThom}
Suppose $P_J$ is homotopic to a contact structure $\xi$ whose Heegaard Floer contact invariant maps nontrivially, $F_{{W},\kappa_J}^*(c(\xi))\neq 0$, and $(\Sigma,\partial \Sigma)\subset (W,Y)$ is a properly embedded $J$-holomorphic curve with null-homologous boundary, positively transverse to $P_J$ such that $\selflink_{P_J}(\partial\Sigma)=\selflink_\xi(\partial\Sigma)$. Then\[-\chi(\Sigma)\le -\chi(\Sigma')\] for any smoothly embedded $\Sigma'\subset W$ with $\partial \Sigma'$ isotopic to $\partial \Sigma$ and $ [\Sigma',\partial \Sigma']= [\Sigma,\partial \Sigma]\in H_2(W,Y)$. 
\end{thm}

\noindent       Theorem \ref{thm:genRelativeSymThom} has a corresponding version for 4-manifolds $W$ equipped with a non-degenerate (not necessarily closed) 2-form $\omega$ and oriented surfaces $\Sigma$ on which $\omega$ is a volume form. In this case, after choosing a metric compatible with $\omega$, we can define a plane field $P_\omega=\ker(\star_Y(\omega|_Y))$. As above, we get a map on Floer cohomology \[F_{W,\kappa_J}^*:\HFa^*(Y)\to \HFa^*(Y_{in})\] where $\kappa_J$ is the $\SpinC$ structure associated to an almost complex structure $J$ compatible with $\omega$ and the metric.

\begin{thm}\label{thm:2-formversion}
    Let $W$ be a 4-dimensional cobordism to a 3-manifold $Y$, equipped with a nondegenerate 2-form $\omega$. Suppose $P_{\omega}$ is homotopic to a contact structure $\xi$ such that $F_{W,\kappa_J}(c(\xi))\neq 0$ and  $(\Sigma,\partial\Sigma)\subset (W, Y)$ is properly embedded surface with null-homologous boundary, positively transverse  to $P_\omega$ for which $\omega|_\Sigma>0$ and $\selflink_{P_\omega}(\partial\Sigma)=\selflink_{\xi}(\partial\Sigma)$. Then \[-\chi(\Sigma)\leq -\chi(\Sigma')\] for any smoothly embedded $\Sigma'\subset W$ with $\partial \Sigma'$ isotopic to $\partial \Sigma$ and $ [\Sigma',\partial \Sigma']= [\Sigma,\partial \Sigma]\in H_2(W,Y)$. 
\end{thm} 
\noindent See Remark \ref{rem:non-degenerate} below for more details.

Applying Theorem \ref{thm:genRelativeSymThom} to a strong symplectic filling, we obtain the relative symplectic Thom conjecture as an easy corollary: 

\begin{corollary}[Relative Symplectic  Thom Conjecture]\label{cor:relSymThom} If $(\Sigma, \partial \Sigma)\subset (W,Y)$ is a smoothly embedded symplectic surface whose boundary is a null-homologous positive transverse link in $(Y,\xi)$, 
 strongly filled by $(W,\omega)$, then \[-\chi(\Sigma)\le -\chi(\Sigma')\]
 for any smoothly embedded $\Sigma'\subset W$ with $\partial \Sigma'$ isotopic to $\partial \Sigma$ and $[\Sigma',\partial \Sigma']= [\Sigma,\partial \Sigma]\in H_2(W,Y)$. 
\end{corollary} 
\noindent Corollary \ref{cor:relSymThom} follows immediately from Theorem \ref{thm:genRelativeSymThom} by observing that since $(W,\omega)$ is a strong filling, the map $F_{{W\setminus B^4},\kappa_\omega}^*:\HFa^*(Y)\to\HFa^*(S^3)$ sends $c(\xi)$ to $c(\xi_{std})$, the non-zero contact invariant of the standard tight contact structure on $S^3$ (see Theorem 2.13 and Remark 2.14 of \cite{Ghiggini-fillability}.) Here the symplectic structure determines an isotopy class of compatible almost complex structures and hence gives a well defined $\SpinC$ structure $\kappa_\omega$ whose first Chern class agrees with the Chern class of the complex tangent bundle \cite{HolDiskFour}.  

A proof of Corollary \ref{cor:relSymThom} first appeared in \cite{Gadgil-Kulkarni}, where it was derived from the closed case by showing that $(W,\Sigma)$ can be embedded into a closed symplectic pair $(\overline{W},\overline{\Sigma})$. More recently, Alfieri and Cavallo gave a different proof in the special case of Stein fillings of rational homology spheres \cite{Alfieri-Cavallo}.

In addition to suggesting Theorem \ref{thm:genRelativeSymThom}, our approach to the relative symplectic Thom conjecture has a significant advantage in that it provides obstructions to a given knot type bounding a symplectic surface in a symplectic 4-manifold.  These utilize the contact ``tau" invariant $\tau_{\xi}(Y,L,[S])$ from \cite{tbbounds} (see Definition \ref{def:contacttau} below). We write $\overline{\xi}$ for the conjugate contact structure.

\begin{thm}\label{thm:TauObstruct}
    Let $(W,\omega)$ be a strong symplectic filling of $(Y,\xi)$ and $L\subset Y$ a null-homologous oriented link. If 
    \[2\tau_{\overline\xi}(Y,L,[S])<\langle c_1(\spinct_\omega),\sigma\rangle -\sigma^2 \]
    then $L$ is not isotopic to the oriented boundary of any smoothly embedded symplectic surface $\Sigma$ without sphere components, whose boundary is positively transverse to $\xi$ and has $[\Sigma_S]=\sigma\in H_2(W).$ 
\end{thm}
\begin{remark} A couple of remarks are in order.
    First, for simplicity we state our results for links that are null-homologous in $Y$ but a version of Theorem \ref{thm:genRelativeSymThom} should hold more generally for rationally null-homologous links. See \cite{Baker-Etnyre, Li-Wu, RelativeAdjunction} for the necessary generalizations of the bounds used in the proof of Theorem \ref{thm:genRelativeSymThom}. 

    Second, the reader might wonder why the  conjugate contact structure $\overline{\xi}$ appears in Theorem \ref{thm:TauObstruct}. This follows from the nature of the relative adjunction formula, in which {\em negative} the Chern class of the tangent bundle naturally computes $-\chi(\Sigma)$ of a symplectic surface.  But $-c_1(\kappa_\omega)$ is the Chern class of $\overline{\kappa}_\omega$, whose map on Floer homology is non-trivial on $c(\overline{\xi})$.
\end{remark}

Theorem \ref{thm:genRelativeSymThom} indeed pushes the relative Thom conjecture beyond the case of symplectic fillings.  For instance, it  immediately implies that properly embedded symplectic surfaces in  half  of the symplectization $(Y\times(-\infty,0],d(e^t\alpha))$ of a contact 3-manifold $(Y,\xi=\ker{\alpha})$ with nonvanishing contact invariant maximize Euler characteristic: 

\begin{corollary}\label{cor:Symplectization} Compact symplectic surfaces with fixed positive transverse boundary in the half-symplectization of a contact structure with non-trivial \os\ invariant maximize Euler characteristic in their relative homology class. \end{corollary} 

\noindent This suggests the following conjecture, asserting that tightness of a contact structure is determined by the symplectic geometry of its symplectization:

\begin{conj}[Tightness is symplectically determined]\label{conj:CharTightness} A contact structure is tight if and only if every compact properly embedded symplectic surface with fixed positive transverse boundary in half its symplectization maximizes Euler characteristic amongst all surfaces with   the same relative homology class. \end{conj}

\noindent This should be compared with another conjectural 4-dimensional interpretation of tightness by the ``slice-Bennequin" inequality \cite[Conjecture 1]{4Dtight}. In light of the conjecture, we find it curious that Theorem \ref{thm:genRelativeSymThom} also implies certain symplectic surfaces in the symplectization of an overtwisted contact structure \emph{are} genus minimizing, provided the contact structure is homotopic to one with nonvanishing contact invariant and the boundary curve has a self-linking number allowed in the tight structure; see Example \ref{ex:OTgenusminimizer}.  On the other hand, Example \ref{ex:OTgenusminimizer} also shows that {\em any} link bounds symplectic surfaces in the symplectization of the overtwisted contact structure on the 3-sphere in the same homotopy class as $\xi_{std}$. These will typically not be genus minimizing. It would be interesting to better understand these phenomena.\\

\noindent {\bf Acknowledgements:} MH was partially supported by NSF grant DMS-2104664 and a Simons Fellowship. 
KR was partially supported by NSF grant DMS-2405452, and received Travel Support from the Simons Foundation. Part of this work was conducted while MH was visiting Cambridge University, and while KR was at Max Planck Institute.  We thank both institutions for their hospitality.

\section{Relative adjunction and genus formulae} 
In what follows, we assume our 3-manifolds are smooth, closed, connected, and oriented. Let $P$ be an  oriented plane field on a 3-manifold $Y$. An oriented knot or link is \emph{positively transverse} to $P$ if its orientation agrees with the positive normal direction to the planes and \emph{negatively transverse} if they disagree.

 \subsection{Self-linking numbers}\label{sec:selflinknum} One can assign a self linking number to a   link $L$ transverse to a plane field $P$,   provided we have a trivialization of $P|_L$.  Any Seifert surface $S$ for a  null-homologous transverse link provides such a trivialization, unique up to homotopy.  Indeed, since $S$ is a surface with boundary, one can choose a nowhere zero section $X$ of the restricted bundle $P|_S$. Composing this section with the exponential map along $\partial S$ determines an oriented push-off $L'$ of $L$. The self-linking of $L$ is the linking number of $L$ with $L'$ or, equivalently, the intersection number of $L'$ with $S$: 
\[\selflink_P(L,S):= \lk(L,L') = S\cdot L'.\]
 Alternatively, one may interpret the self-linking number as the evaluation of a relative Chern class. Since $L$ is transverse, its Seifert framing yields a homotopy class of nowhere zero section $s_\nu$ of $P|_L$.  The self-linking number of a transverse link is the obstruction to extending this section over all of $S$ or, equivalently, the Chern class of $P|_S$ \emph{relative} to $s_\nu$ evaluated on the class of the coherently oriented Seifert surface: \[\selflink_P(L,S):=\mp\langle c_1(P|_S,s_\nu),[S, \partial S]\rangle \] 
 where we have the negative sign when $L$ is \emph{positively} transverse, following the sign conventions of \cite{Geiges}. 
Note that reversing the orientation of $L$ changes $L$ from positively to negatively  transverse (and vice versa), and also changes the sign of $S$, so that the self-linking number is independent of the chosen orientation for $L$. 

\begin{figure}
     \centering
     \includegraphics[width=0.5\linewidth]{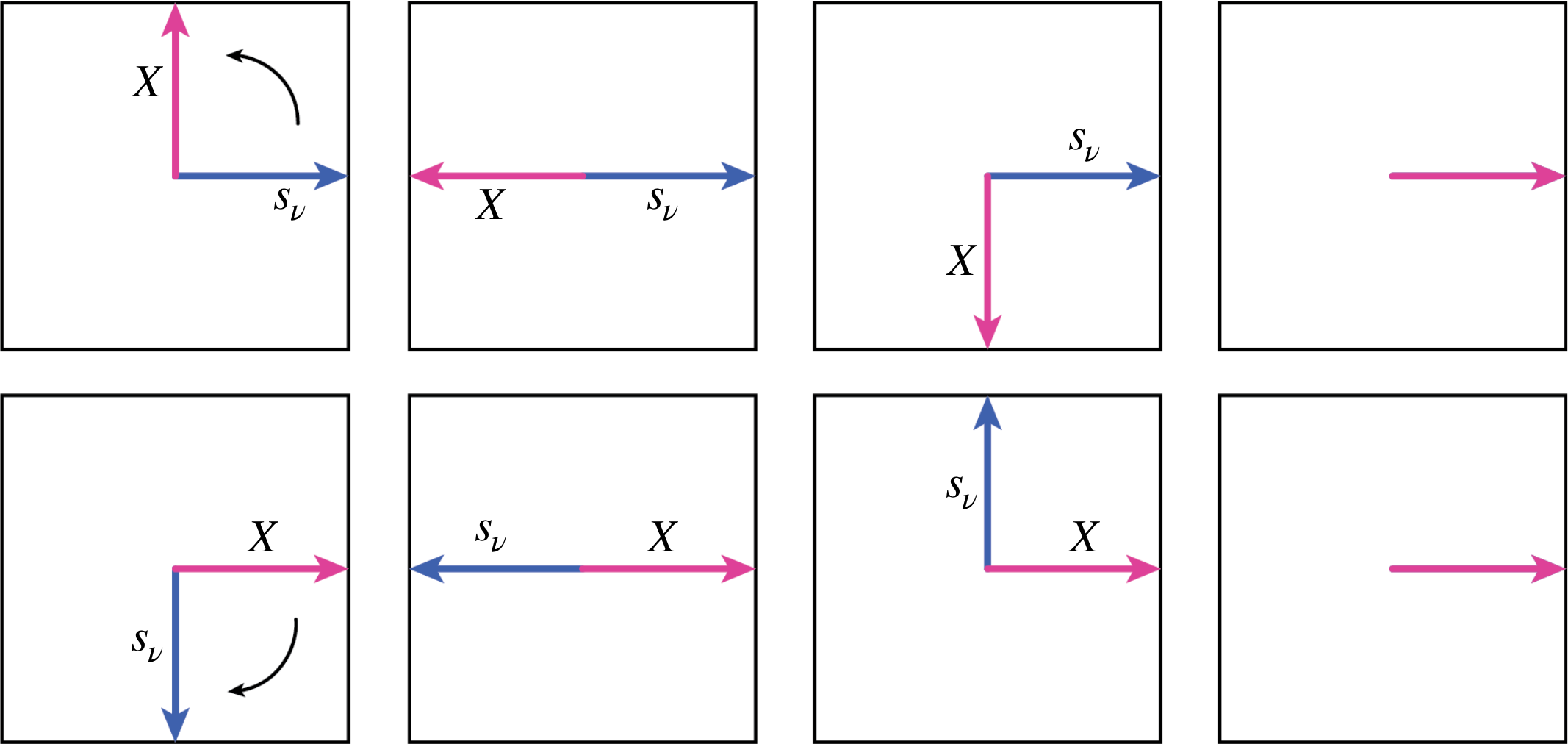}
     \caption{The top sequence shows $X$ winding once counterclockwise with respect to $s_\nu$. The bottom sequence is the same, but viewed as $s_\nu$ winding with respect to $X$ --- from this perspective $s_\nu$ winds once clockwise.}
     \label{fig:winding} 
 \end{figure}

To show these definitions agree, we interpret both as winding numbers. Assume $L$ is positively transverse. The linking number $\lk(L,L')$ is given by the winding of the section $X$, defined above, with respect to the trivialization of $P|_L$ determined by the Seifert framing, $s_\nu$. At the same time, the Chern class evaluation $\langle c_1(P|_S,s_\nu),[S, \partial S]\rangle$ is the obstruction to extending the section $s_\nu$ to a nontrivial section over all of $S$. Thus, if we use $X$ to  trivialize $P|_L$, we can measure the winding of $s_\nu$ with respect to $X$. Since we have reversed the roles of $X$ and $s_\nu$, this winding number has the opposite sign of $\lk(L,L')$, see Figure \ref{fig:winding}.

\subsection{Relative adjunction formulas}
Let $W$ be a smooth oriented 4-manifold with boundary, equipped with an almost complex structure $J: TW\rightarrow TW$. The almost complex structure induces a plane field $P_J$ on any connected component $Y$ of the boundary, given by the intersection of the $3$-dimensional subbundle $TY\subset TW$ with its image under $J$; that is, $P_J:=TY\cap JTY.$

Now let $(\Sigma,\partial\Sigma)\subset (W,Y)$ be a smooth and properly embedded oriented surface  whose boundary is positively transverse to $P_J$.  We consider two complementary sections of $TW|_{\Sigma}$. First, we have $s_\tau$, which is the section of unit tangent vectors to $\partial \Sigma$ determined by the orientation. Second, assuming that $\partial \Sigma$ is null-homologous in $Y$, a Seifert surface $S$ determines a section $s_\nu$  of the plane field $P_J$, as above. 

If $\Sigma$ is $J$-holomorphic, the restriction $TW|_\Sigma$ splits as a sum of complex line bundles $T\Sigma\oplus \nu\Sigma$. The Chern class of $TW|_\Sigma$ \emph{relative} to the trivialization along the boundary given by the pair of sections $\{s_\tau,s_\nu\}$ measures the obstruction to extending this pair of sections to linearly independent sections over all of $\Sigma$.  We have the following relative version of the well-known adjunction formula from algebraic geometry.  

\begin{prop}[The Relative Adjunction Formula] \label{prop:adjformula} Let $(W,J)$ be an almost complex 4-manifold with boundary, and let $(\Sigma,\partial\Sigma)\subset (W,Y)$ be a properly embedded $J$-holomorphic curve whose boundary is transverse to the plane field $P_J$. {For any non-zero section $s$ of   $\nu\Sigma$ along $\partial \Sigma$}, we have:
\[ c_1(TW|_\Sigma, s_{\tau}\oplus s) = c_1(T\Sigma,s_{\tau})+ c_1(\nu\Sigma,s).\]
\end{prop}
\begin{proof}  Like the adjunction inequality for closed surfaces, the relative adjunction inequality is a consequence of the Whitney product formula equating the total Chern class of a direct sum of bundles with the product of the total Chern classes of the summands. We employ an analogous  formula for the relative Chern classes. One can geometrically derive the formula by appealing to the obstruction theoretic interpretation of the classes involved; alternatively, one can use the general theory of relative characteristic classes due to Kervaire \cite{Kervaire}.  The appendix of that work, in particular, addresses relative Chern classes, with the Whitney product formula (referred to as ``Whitney duality") proved on page 549.  The above formula is then an easy special case. 
\end{proof}

\noindent Note that a somewhat different relative adjunction formula figures prominently in Hutching's calculation of the ECH index, see \cite[Proposition 4.9]{HutchingsRevisit} for a statement and \cite[Section 3]{Hutchings} for its proof.

Evaluating both sides of the relative adjunction formula on the fundamental class $[\Sigma,\partial\Sigma]$, we recover a formula that computes the Euler characteristic of $\Sigma$, analogous to the genus formula for a closed complex curve.  With the aid of the self-linking number, we can manipulate this into a more useful genus formula involving only absolute homology and cohomology classes:

\begin{prop}[Genus Formula]\label{prop:GenusFormula} 
Let $(W,J)$ be an almost complex 4-manifold with boundary, and let $(\Sigma,\partial\Sigma)\subset (W,Y)$ be a properly embedded $J$-holomorphic curve whose boundary is positively transverse to the plane field $P_J$.  If $\partial\Sigma$ is nullhomologous in $Y$, then
\[ -\langle c_1(W,J),[\Sigma_{S}]\rangle +  [\Sigma_{S}]\cdot[\Sigma_S] + \selflink_{P_J}(\partial\Sigma,S) = -\chi(\Sigma),\]
where $S$ denotes a Seifert surface for $\partial\Sigma$ used to compute the self-linking, and $\Sigma_{S}$ is the closed surface obtained as the union $\Sigma\cup -S$.
\end{prop} 
\noindent Before proving the proposition, we note that both the self-linking number and the homology  class $[\Sigma_S]$ depend on the relative homology class  $[S,\partial S]$, but the left hand side of the formula is independent of this choice because $c_1(W,J)|_{Y}=c_1(P_J)$.

\begin{proof} Let $s_\tau$ be the unit tangent vectors to $\partial\Sigma$ and let $s_\nu$ be the Seifert framing as in Section \ref{sec:selflinknum}.
    We use Proposition \ref{prop:adjformula} to evaluate $c_1(TW|_\Sigma, s_\tau\oplus s_\nu)$ on $[\Sigma,\partial\Sigma]$. This gives 
    \begin{equation}\label{eqn:relChernEval}
       \langle c_1(TW|_\Sigma, s_\tau\oplus s_\nu),[\Sigma,\partial\Sigma]\rangle =\chi(\Sigma)+ \langle c_1(\nu\Sigma,s_\nu),[\Sigma,\partial\Sigma] \rangle,
    \end{equation}
\noindent since evaluating  the relative Chern class $c_1(T\Sigma,s_\tau)$ on $[\Sigma,\partial\Sigma]$ recovers the Euler  characteristic of the surface. 

Next, we consider the rightmost term in Equation \ref{eqn:relChernEval}. We claim this relative Chern class evaluation recovers the self-intersection number of $\Sigma_S$: 
\begin{equation}\label{eqn:Selfint}
    \langle c_1(\nu\Sigma,s_\nu),[\Sigma,\partial \Sigma]\rangle = [\Sigma_S]\cdot[\Sigma_S].
\end{equation}
To see this, observe that $\langle c_1(\nu\Sigma,s_\nu),[\Sigma,\partial \Sigma]\rangle$ can be computed as the intersection number of $\Sigma$ with a parallel pushoff in a direction normal to  $\Sigma$ that agrees with $s_\nu$ on the boundary. At the same time, we claim that the Seifert surface doesn't contribute anything to the self-intersection number of $\Sigma_S$. To see this, extend $s_\nu$ over $S$ and form a pushoff $S'$ in the direction of $s_\nu$ for $S$. Since $s_\nu$ agrees with the Seifert framing for $S$
, we have $S\cdot S'=\lk(\partial S,\partial S')=0$. 
Thus, $\Sigma\cdot \Sigma'=  \Sigma_S\cdot (\Sigma_S)'$ where $(\Sigma_S)'=\Sigma'\cup -S'$. Equation \ref{eqn:Selfint} now follows. 

Putting Equations \ref{eqn:relChernEval} and \ref{eqn:Selfint} together we now have:
\begin{equation}\label{eqn:RHS}
     \langle c_1(TW|_\Sigma, s_\tau\oplus s_\nu),[\Sigma,\partial\Sigma]\rangle =\chi(\Sigma)+[\Sigma_S]\cdot[\Sigma_S]
\end{equation}

We turn our focus now to the left hand side of Equation \ref{eqn:RHS}. Here, we use a key naturality property of relative characteristic classes, which states that their image under the map from relative cohomology to absolute cohomology is the (non-relative) characteristic class \cite[Lemma 11.4]{Kervaire}. For our purposes, this means $c_1(W,J)=j^*c_1(TW,s_\tau\oplus s_\nu)$ where $j^*:H^2(W,\partial\Sigma)\to H^2(W)$ is the map in the long exact sequence of the pair $(W,\partial\Sigma)$. From this we compute: 
\begin{align*}
\langle c_1(W,J),[\Sigma_{S}]\rangle 
& =\langle j^*c_1(TW, s_\tau\oplus s_\nu), [\Sigma_S]\rangle\\
&=\langle c_1(TW, s_\tau\oplus s_\nu), j_*[\Sigma_S]\rangle\\
&=\langle c_1(TW, s_\tau\oplus s_\nu), [\Sigma,\partial\Sigma]-[S,\partial S]\rangle\\
&=\langle c_1(TW, s_\tau\oplus s_\nu), [\Sigma,\partial\Sigma]\rangle -\langle c_1(TW, s_\tau\oplus s_\nu), [S,\partial S]\rangle
\end{align*}
where the final lines follow by observing that the map \[j_*:H_2(W)\to H_2(W,\partial\Sigma)\] sends $[\Sigma_S]$ to the difference of relative classes: $[\Sigma,\partial\Sigma]-[S,\partial S]$. 

Rearranging the computation gives:
\begin{equation}\label{eqn:LHS}
  \langle c_1(TW, s_\tau\oplus s_\nu), [\Sigma,\partial\Sigma]\rangle =\langle c_1(W,J),[\Sigma_{S}]\rangle + \langle c_1(TW, s_\tau\oplus s_\nu), [S,\partial S]\rangle. 
\end{equation}

 To complete the argument, we show the rightmost term in Equation \ref{eqn:LHS} is \emph{minus} the self-linking number. To see this, it is helpful to restrict $TW$ to $Y$ before restricting to $S$.  This is because the restriction $TW|_Y$ splits as $P_J^\perp\oplus P_J$ where $P_J^\perp$ denotes the complement in $TW|_Y$. Note that the complex line bundle $P_J^\perp$ has a non-zero section $s_{out}$ defined by the outward normal vectors to $Y$.  
 
 We would like to apply the relative Whitney sum formula to $P_J^\perp\oplus P_J$.  By definition, $s_\nu$ is a nowhere zero section of $P_J$ along $\partial S$. The section $s_\tau$, on the other hand, may not live in $P_J^\perp$, so the Whitney sum formula can't be directly applied. Observe, however, that $s_\tau$ is homotopic to $s_{out}$ along $\partial S$ via the homotopy $s_t=(1-t)s_\tau+ ts_{out}$.  As $s_\tau$ is  linearly independent from $s_{out}$, this homotopy is through non-vanishing sections. Since relative Chern classes are invariant under homotopy of framings \cite[pgs 520-521]{Kervaire}, we have the following computation:
\begin{align*}
\langle c_1(TW, s_\tau\oplus s_\nu),(i_W)_*[S,\partial S]\rangle &= \langle c_1(TW|_Y,s_\tau\oplus s_\nu),(i_Y)_*[S,\partial S]\rangle \\
&= \langle c_1(TW|_Y,s_{out}\oplus s_\nu),(i_Y)_*[S,\partial S]\rangle \\
&= \langle c_1(P_J^\perp\oplus P_J ,s_{out}\oplus s_\nu),[S,\partial S]\rangle\\
&=\langle c_1(P_J^\perp,s_{out}),[S,\partial S]\rangle+\langle c_1(P_J,s_\nu),[S,\partial S]\rangle\\
&=\langle c_1(P_J^\perp,s_{out}),[S,\partial S]\rangle -\selflink_{P_J}(\partial\Sigma, S)\\
&=-\selflink_{P_J}(\partial\Sigma, S).
\end{align*}

\noindent Here the final line follows because the relative Chern class $c_1(P_J^\perp,s_{out})$ is zero, since the section $s_{out}$ along $\partial S$ extends to a non-vanishing section along all of $Y$. Using this computation to replace the rightmost term in Equation \ref{eqn:LHS}, we obtain:  
\begin{equation}\label{eqn:absChernclass}
\langle c_1(TW,s_\tau\oplus s_\nu),[\Sigma,\partial\Sigma]\rangle=\langle c_1(W,J),[\Sigma_S]\rangle - \selflink_{P_J}(\partial\Sigma,S). 
\end{equation}
Finally, we plug Equation \eqref{eqn:absChernclass} into the left hand side of Equation \eqref{eqn:RHS} to conclude: \[\langle c_1(W,J),[\Sigma_S]\rangle - \selflink_{P_J}(\partial\Sigma,S) =\chi(\Sigma)+[\Sigma_S]\cdot[\Sigma_S].\]
\end{proof}

\begin{remark}\label{rem:non-degenerate}
    Note that the results above  only used the fact that the tangent bundle to $W$ splits along $\Sigma$ as a sum of complex line bundles with preferred boundary sections.  If $\Sigma\subset(W,\omega)$ is surface in a $4$-manifold equipped with a non-degenerate (not necessarily symplectic) 2-form, and $\omega|_{\Sigma}>0$,  
 then the symplectic vector bundle $TW|_{\Sigma}$ splits as a sum of symplectic vector bundles, and hence as a sum of complex line bundles, see e.g. \cite[Theorem 2.1]{McCarthyWolfson}.  Therefore we have an analogous version of both Proposition \ref{prop:adjformula} and \ref{prop:GenusFormula}, where $\Sigma$ is assumed to be symplectic with respect to a non-degenerate 2-form on $W$, and with boundary transverse to the plane field on $Y$ given by the kernel of $\star_Y(\omega|_Y)$. 
 \end{remark}

\section{The relative almost complex Thom conjecture}\label{sec:relThom}

  In this section, we prove Theorem \ref{thm:genRelativeSymThom} and its corollaries. To this end, we first briefly recall the definition of the relevant $\tau$ invariants. We refer the reader to \cite{tbbounds, 4Dtight, RelativeAdjunction} for further details.

In \cite{RelativeAdjunction} the authors associate a collection of Heegaard Floer theoretic invariants to a  null-homologous knot $K$ in a 3-manifold $Y$. A choice of  Seifert surface for $K$ determines an Alexander filtration $\mathcal{F}_m(Y,K,[S])$ of the Heegaard Floer chain complex $\CFa(Y)$. The values of the Alexander filtration depend only on the relative homology class of the chosen surface $S$. For each $m\in\Z$ we have an inclusion $\iota_m:\mathcal{F}_m(Y,K,[S])\hookrightarrow \CFa(Y)$ which induces a map $I_m$ on homology. 

% \begin{defn} For $\alpha\in\HFa(Y)$ we define
% \begin{align*}
%     \tau_{\alpha}(Y,K,[S])&:=\min\{r\in\Q\mid \alpha\in\image(I_r)\}.
% \end{align*}
% \end{defn}
\begin{defn}\label{def:contacttau} For $\alpha\in\HFa(Y)$ and $\varphi\in\HFa^*(Y)\cong \HFa(-Y)$ we define
\begin{align*}
    \tau_{\alpha}(Y,K,[S])&:=\min\{m\in\Z\mid \alpha\in\image(I_m)\}\\
    \tau^*_{\varphi}(Y,K,[S])&:=\min\{m\in\Z\mid \exists \beta\in \image(I_m) \text{ such that } \langle \varphi,\beta\rangle \neq 0\}
\end{align*}
where $\langle -,-\rangle$ denotes the duality pairing between $\HFa^*(Y)$ and $\HFa(Y)=\HFa_*(Y)$.
\end{defn}
\noindent The authors  extend these invariants to links $L\subset Y$ by ``knotifying"  $L$ to obtain a knot $\kappa(L)$ in $Y\#^{|L|-1}S^1\times S^2$. For each Floer class $\Theta\in \HFa(\#^{|L|-1}S^1\times S^2)$ we have: \[\tau_{\alpha\otimes\Theta}(Y,L,[S]):=\tau_{\alpha\otimes\Theta}(Y\#^{|L|-1}S^1\times S^2,\kappa(L),[S_{\kappa}])\] where $S_\kappa$ is the Seifert surface for $\kappa(L)$ obtained from $S$ by attaching bands which pass through the $1$-handles used in the knotification procedure.  

In \cite[Definition 1.4]{4Dtight} the authors apply these constructions to the contact class $c(\xi)$ which is naturally a Floer cohomology class in $\HFa^*(Y)$. We define: 
\begin{equation}\label{eqn:tauxi}
    \tau_{\xi}(Y,L,[S]):=\tau^*_{c(\xi)\otimes c(\xi_{std})}(Y\#^{|L|-1}S^1\times S^2,\kappa(L),[S_{\kappa}])
\end{equation}
where $\xi_{std}$ is the unique tight contact structure on $\#^{|L|-1}S^1\times S^2.$ 

Let $L^r$ denote $L$ with its orientation reversed. We will need the following symmetry property of $\tau_{\xi}(Y,L^r,[-S])$:

\begin{prop}\label{prop:OrientRev}
    Let $L$ be an oriented link in a contact 3-manifold $(Y,\xi)$. Then, \[\tau_{\xi}(Y,L^r,[-S])=\tau_{\overline{\xi}}(Y,L,[S])\]
    where $\overline{\xi}$ is the conjugate contact structure on $Y$ i.e. the same two plane field, but with reversed orientation. 
\end{prop}
\begin{proof}
    This follows from \cite[Proposition 3.6]{Raoux} which shows  \[\CFKinf(Y,K,\spinc)\cong \CFKinf(Y,K^r,\overline{\spinc}).\] Applying this to the knotification of $L$ we have \[\CFK_{\infty}^*(Y\#^{|L|-1} S^1\times S^2,\kappa(L), \spinc_{\xi\#\xi_{std}})\cong \CFK_{\infty}^*(Y\#^{|L|-1} S^1\times S^2,\kappa(L)^r, \overline{\spinc_{\xi\#\xi_{std}}}).\] Note that $\overline{\spinc_{\xi\#\xi_{std}}}=\spinc_{\overline{\xi}\#\xi_{std}}$ since $\xi_{std}$ on $\#^{|L|-1} S_1\times S_2$ is self-conjugate. Furthermore, under this identification of doubly-filtered complexes, generators representing $c(\xi\#\xi_{std})$ are sent to generators representing $c(\overline{\xi}\#\xi_{std})$.
\end{proof}

We now collect the two key properties of $\tau_\xi$ used in our proof Theorem \ref{thm:genRelativeSymThom}.  The first is a Bennequin-type bound:

\begin{theorem}[$\tau$-Bennequin Inequality] Let $(Y,\xi)$ be a contact 3-manifold with non-vanishing \os\ invariant $c(\xi)\ne0\in \HFa^*(Y)$.  If $\mathcal{L}$ is a Legendrian link whose underlying smooth link type $L$ has $|L|$ components then 
\begin{equation}\label{eqn:tauBennequin}
\tb_\xi(\mathcal{L})+\rot_{\xi}(\mathcal{L}, S)\leq 2\tau_{\xi}(Y,L,[S])-|L|.
\end{equation}
If $L$ is a positively transverse link, then \begin{equation}\label{eqn:selflink}
\selflink_\xi(L,S)\le 2\tau_{\overline{\xi}}(Y,L,[S])-|L|
\end{equation}
\end{theorem}

\begin{proof}
 Equation \ref{eqn:tauBennequin} is established for knots in \cite{tbbounds} and for links in \cite{4Dtight}.

 For the second part, note that a positively transverse link $L$ can be approximated by a Legendrian $\mathcal{L}$ whose positive transverse push-off is $L$ so that \[\selflink_\xi(L,S)=\tb_{\xi}(\mathcal{L})-\rot_{\xi}(\mathcal{L},S)=\tb_{\xi}(\mathcal{L}^r)+\rot_{\xi}(\mathcal{L}^r,-S)\] where the second equality follows by reversing the orientation of $L$ \cite[pgs. 130-131]{Etnyre-survey}. It follows from Equation \eqref{eqn:tauBennequin} and Proposition \ref{prop:OrientRev} that \[\selflink_\xi(L,S)\leq 2\tau_{\xi}(Y,L^r,[-S])-|L|=2\tau_{\overline{\xi}}(Y,L,[S])-|L|.\]
\end{proof}

The other property of $\tau_\xi(Y,L,[S])$ we will need is that it bounds the genera of smoothly embedded surfaces in cobordisms that map the contact class non-trivially:
\begin{theorem} [Relative Adjunction Inequality] \label{prop:linkadjunction}
    Let $W$ be a  4-dimensional cobordism between 3-manifolds, $Y_1$ and $Y_2$ such that $F_{{W}, \spinct}^*:\HFa^*(Y_2)\to\HFa^*(Y_1)$ satisfies  $F_{{W}, \spinct}^*(c(\xi))\ne 0$.
    If $(\Sigma,\partial\Sigma)\subset (W,Y_2)$ is a smooth, properly embedded surface whose boundary has $|\partial\Sigma|$ components, we have: 
    \begin{equation}\label{eqn:relative adjunction}
         \langle c_1(\spinct),[\Sigma_S]\rangle +[\Sigma_S]\cdot[\Sigma_S]+2\tau_{\xi}(Y_2,\partial\Sigma,[S])-|\partial\Sigma| \leq  -\chi(\Sigma).
     \end{equation}
   
\end{theorem}

\begin{proof}
    If $L$ is a knot, this is a special case of \cite[Theorem 4.2]{RelativeAdjunction}. The general case essentially follows the proof of \cite[Theorem 5.17]{RelativeAdjunction}, c.f. \cite[Section 4.1]{Rational4genus}.  
    \end{proof}

We are now ready to prove Theorem \ref{thm:genRelativeSymThom}.
\begin{proof}[Proof of Theorem {\ref{thm:genRelativeSymThom}}]
 Let $(W,J)$ be an almost complex 4-manifold with boundary $\partial W=-Y_{in}\sqcup Y$, regarded as a cobordism from $Y_{in}$ to $Y$.  Suppose $\Sigma\subset W$ is a properly embedded $J$-holomorphic curve whose boundary $\partial \Sigma$ lies in $Y$, where it is positively transverse to the induced plane field $P_J=TY\cap JTY$. According to Proposition \ref{prop:GenusFormula}, the Euler characteristic of $\Sigma$ is given as
$-\chi(\Sigma)$:  
\begin{align*}
    -\chi(\Sigma)&=-\langle c_1(W,J),[\Sigma_S]\rangle +[\Sigma_S]^2 +\selflink_{P_J}(\partial\Sigma,S)\\
    &=-\langle c_1(W,J),[\Sigma_S]\rangle +[\Sigma_S]^2 +\selflink_{\xi}(\partial\Sigma,S)
\end{align*}
\noindent since, by assumption, $\selflink_{P_J}(\partial \Sigma,S)=\selflink_\xi(\partial \Sigma,S)$.

 Since the contact class of $\xi$ is non-zero, we can apply the $\tau$-Bennequin inequality \eqref{eqn:selflink} to $\partial\Sigma$: 
 \[ \selflink_{\xi}(\partial\Sigma,S)\le 2\tau_{\overline{\xi}}(Y,\partial \Sigma,[S])-|\partial \Sigma|.\] This implies \[-\chi(\Sigma)\leq -\langle c_1(W,J),[\Sigma_S]\rangle +[\Sigma_S]^2 +2\tau_{\overline{\xi}}(Y,\partial \Sigma,[S])-|\partial \Sigma|\] 
On the other hand, since $F_{W,\kappa_J}^* (c(\xi))\ne0 $ we also have $F_{W,\overline{\kappa}_J}^* (c(\overline{\xi}))\ne0 $. The relative adjunction inequality \eqref{eqn:relative adjunction} applied to any surface $(\Sigma',\partial\Sigma')$ with $\partial\Sigma'$ isotopic to $\partial\Sigma$ and representing the same relative homology class as $[\Sigma,\partial\Sigma]$, satisfies
 \[ \langle c_1(\overline{\kappa}_J),[\Sigma'_{S}]\rangle +[\Sigma'_{S}]\cdot[\Sigma'_{S}]+2\tau_{\overline{\xi}}(Y,\partial\Sigma',[S])-|\partial\Sigma'|\leq -\chi(\Sigma').  \]

\noindent Now observe that $c_1(\overline{\kappa}_J)=-c_1(W,J)$. We can therefore string the inequalities above together to yield:
\[-\chi(\Sigma)\le -\chi(\Sigma'),\]
for any smooth surface $\Sigma'$ in the same relative homology class as $\Sigma$ with  $\partial \Sigma'\simeq \partial \Sigma$. The result now follows.
\end{proof}

Alternatively, we can consider a surface which is symplectic for a (not necessarily closed) non-degenerate 2-form $\omega\in \Omega^2(W)$, according to Remark \ref{rem:non-degenerate}. 

\begin{proof}[Proof of Theorem \ref{thm:2-formversion}]
    Theorem \ref{thm:2-formversion} follows the same argument as Theorem \ref{thm:genRelativeSymThom}, once the adjunction formula is established. See Remark \ref{rem:non-degenerate}. 
\end{proof}

 Given a smooth null-homologous link $L\subset Y$ in a strong filling $(W,\omega)$  of $(Y,\xi)$ we also obtain an obstruction to $L$ bounding a symplectic surface:

\begin{proof}[Proof of Theorem \ref{thm:TauObstruct}] Suppose $(W,\omega)$ is a strong filling of $(Y,\xi)$ and $L$ is isotopic to the oriented positive transverse  boundary of a properly embedded symplectic surface $\Sigma$ in $(W,\omega)$. Let $S$ be a Seifert surface such that $\sigma=[\Sigma_S].$

For each connected component $\Sigma_i$ of $\Sigma$ we have \[|\partial\Sigma_i|-\chi(\Sigma_i)=2g(\Sigma_i)+2|\partial\Sigma_i|-2\geq 0\] provided $\Sigma$ has no sphere components. Summing up over connected components implies $|\partial\Sigma|-\chi(\Sigma)\geq 0$.

Combining this estimate with Theorem \ref{thm:genRelativeSymThom} now implies: \[-\langle c_1(\kappa_\omega),[\Sigma_S]\rangle +[\Sigma_S]^2+2\tau_{\overline{\xi}}(Y,L,[S])=|\partial\Sigma|-\chi(\Sigma)\geq 0.\] 
Thus, a link with $2\tau_{\overline{\xi}}(Y,L,[S])< \langle c_1(\kappa_\omega),\sigma\rangle -\sigma^2$ violates the above inequality. 
\end{proof}

For example, if $(Y,\xi)$ is a contact rational homology sphere with $c(\xi)\neq 0$ and $K\subset Y$ has $\tau_{\overline{\xi}}(Y,K)<0$ then $K$ cannot bound a symplectic surface in any rational homology ball filling of $(Y,\xi)$. Furthermore, since $\tau_{\overline{\xi}}(Y,K)\leq \max_{\alpha\in\HFa(Y)} \tau_{\alpha}(Y,K)$ a knot $K\subset Y$ whose maximal $\tau$ invariant is negative cannot bound a symplectic surface in any rational homology ball filling of any contact structure on $Y$. Examples are easy to construct. Start with any knot $K\subset Y$ and connect sum with a negative trefoil $-T_{2,3}\subset S^3$. Additivity of the $\tau$ invariant implies $\tau_{\alpha}(Y\# S^3, K\# -T_{2,3})= \tau_{\alpha}(Y,K)-1$ for all $\alpha\in\HFa(Y)$, so all the invariants shift down by 1. Thus, connect summing sufficiently many negative trefoils we obtain a knot that doesn't bound a symplectic surface.

Now we turn to Corollary \ref{cor:Symplectization}. 
The \emph{symplectization} of a contact manifold $(Y,\xi)$ is the symplectic manifold $(Y\times \R, d(e^t\alpha))$, where $\alpha$ is a 1-form on $Y$ such that $\ker(\alpha)=\xi$. For any choice of closed interval $[a,b]\subset \R$ the compact submanifold $Y\times [a,b]\subset Y\times \R$ has the structure of a Stein cobordism from $(Y\times \{a\},\xi)$ to $(Y\times\{b\},\xi)$.  We'll refer to $Y\times (-\infty,0]$ with the symplectic structure restricted from that above as the {\em half-symplectization} of $(Y,\xi)$. 

\begin{prop}\label{prop:half-symp} Let 
     $Y\times (-\infty,0]$ be the half-symplectization of  a contact 3-manifold $(Y,\xi)$ with nontrivial contact invariant. If $\Sigma$ is a compact, properly embedded surface with positive transverse boundary in $Y\times\{0\}$ then $\Sigma$ maximizes Euler characteristic in its relative homology class. 
\end{prop}

\begin{proof} If $\Sigma$ is a compact surface properly embedded in $Y\times (-\infty,0]$ with positive transverse boundary in $Y\times \{0\}$. Then for any $a<0$ such that $\Sigma\subset Y\times[a,0]$ we show $\Sigma$ is genus minimizing in its relative homology class.
 Since the product manifold $Y\times [a,0]$ maps every Floer class by the identity, the non-zero contact class $c(\xi)\in\HFa(-Y,\spinc_\xi)$ maps nontrivially. Theorem \ref{thm:genRelativeSymThom} now implies that it has maximal Euler characteristic among all surfaces in the same  relative homology class in $Y\times [a,0]$ with the given boundary.  
\end{proof}
\noindent Proposition \ref{prop:half-symp} immediately implies Corollary \ref{cor:Symplectization}.

From here we can fairly easily prove one direction of Conjecture \ref{conj:CharTightness}:

\begin{prop} If every properly embedded symplectic surface in the symplectization of a contact manifold maximizes Euler characteristic in its relative homology class then the contact structure is tight.
\end{prop}
\begin{proof}  We prove the contrapositive.  Suppose then, that a contact structure is overtwisted.  We will find a properly embedded symplectic surface in its symplectization which does not maximize Euler characteristic for its boundary. To do this, we observe that the overtwisted contact structure on the $3$-sphere with trivial Hopf invariant has an open book decomposition with connected binding and whose page is not $4$-genus minimizing (e.g. $T_{p,q}\#-T_{p,q}\#K$, where $K$ is a fibered knot with Hopf invariant opposite that of $-T_{p,q}$).  Now we appeal to Proposition \ref{prop:symplecticpage} (below) which states that subsurfaces of pages of open books are properly isotopic to symplectic surfaces in the symplectization. \end{proof}

\begin{prop}\label{prop:symplecticpage} 
    Subsurfaces of pages of open books are properly isotopic to symplectic surfaces in the symplectization.
\end{prop} 
\begin{proof}
    A subsurface of a page of an open book is  isotopic to a Legendrian ribbon in the contact structure associated to the open book. In \cite[Example 4.7]{HaydenStein}, Hayden shows how to push such a ribbon into the symplectization to be ascending which, by \cite[Lemma 5.1]{HaydenStein}, suffices to isotope the surface to be symplectic.  
\end{proof}

On the other hand, in the symplectization of certain overtwisted contact manifolds, some symplectic surfaces \emph{are} genus minimizers. Specifically, if the contact structure is \emph{homotopic} to one with nonvanishing contact invariant, Theorem \ref{thm:genRelativeSymThom} shows that a symplectic surface maximizes Euler characteristic whenever the self-linking numbers of the boundary curve with respect to both plane fields agree.

\begin{example}\label{ex:OTgenusminimizer}
Let $(Y,\xi_{OT})$ be a contact 3-manifold with overtwisted contact structure. Suppose that $\xi_{OT}$ is homotopic to a contact structure $\eta$ with $c(\eta)\neq 0$. If $\Sigma$ is a properly embedded symplectic surface with positively transverse boundary and $\selflink_{P_J}(\partial\Sigma)=\selflink_{\xi}(\partial\Sigma)$ then Theorem \ref{thm:genRelativeSymThom} implies that $\Sigma$ maximizes  Euler characteristic.
\medskip

\noindent An abundant and illuminating source of  this phenomenon  can be found in the 3-sphere.  If $\xi_{0}$ is the overtwisted contact structure on $S^3$ with the same Hopf invariant as the tight contact structure $\xi_{std}$, one can show that {\em any} Seifert surface for {\em any} link is isotopic onto the page of an open book supporting $\xi_0$ \cite{Lyon}.  It follows from Proposition \ref{prop:symplecticpage} that any Seifert surface can therefore be properly isotoped to be symplectic. Such surfaces will typically \emph{not} maximize Euler characteristic, even in the 3-sphere itself. They will, however, maximize Euler characteristic when their self-linking numbers can be realized in the standard tight contact structure, since those self-linking numbers are constrained by the  Bennequin inequality.  In either case, the self-linking number is given as $-\chi(\Sigma)$ of the Seifert surface in question.   Thus, when the Bennequin bound is sharp, we can produce an Euler characteristic maximizing symplectic surface in the half-symplectization of $\xi_0$.  In an upcoming paper, we establish a type of h-principle for transverse knots which allows us to extend such a symplectic surface to one transverse to the tight contact structure on the 3-sphere.  It would be very interesting if this construction could shed light on the conjectured equivalence between strong quasipoisitivity and sharpness of the Bennequin bound \cite{K3List}.

\end{example}

\bibliographystyle{plain}
\bibliography{mybib}

\end{document}